\documentclass[pmlr]{jmlr}

\usepackage{longtable}

\usepackage{array}
\usepackage{amscd}
\usepackage{comment}
\usepackage{stmaryrd}
\usepackage{xcolor}
\usepackage{enumerate}
\usepackage{array}
\usepackage{thmtools}

\newcommand{\X}{{\cal X}}

\def\S{{\cal S}}

\newcommand{\Z}{{\cal Z}}

\newcommand{\argmax}{\mathop{\rm arg\, max}}

\newcommand{\supp}{\mathop{\rm supp}}

\newcommand{\Rank}{\mathop{\rm Rank}}

\newcommand{\RR}{\mathbb{R}}
\newcommand{\NN}{\mathbb{N}}

\newcommand{\roc}{\rm ROC}

\newcommand{\mv}{\rm MV}

\newcommand{\by}{\mathbf{y}}

\newcommand{\bZ}{\mathbf{Z}}

\newcommand{\iid}{\textit{i.i.d.}}
\newcommand{\cdf}{\textit{c.d.f.}\;}

\newcommand{\rv}{\textit{r.v.}}
\newcommand{\ie}{\textit{i.e.}\;}

\newcommand{\wrt}{\textit{w.r.t.}\;}

\newcommand{\eg}{\textit{e.g.}}

\newcommand{\resp}{\textit{resp.}}

\author{\Name{Myrto Limnios} \Email{ myrto.limnios@ens-paris-saclay.fr} \AND
		\addr Universit\'e Paris-Saclay, ENS Paris-Saclay\\
		 \addr CNRS UMR 9010, Centre Borelli, 91190 Gif-sur-Yvette, France
		\AND
		 \Name{Nathan Noiry}  \Email{nathan.noiry@telecom-paris.fr} \AND
		 \Name{Stephan Cl\'emen\c{c}on}  \Email{stephan.clemencon@telecom-paris.fr} \AND
		\addr Telecom Paris, LTCI, Institut Polytechnique de Paris\\
	 \addr19 place Marguerite Perey, Palaiseau, 91120, France}

\editors{Nuno Moniz, Paula Branco, Lu\'{i}s Torgo, Nathalie Japkowicz, Michał Woźniak and Shuo Wang.}

\title[Learning to Rank Anomalies]{Learning to Rank Anomalies: Scalar Performance Criteria and Maximization of Two-Sample Rank Statistics }

\jmlrvolume{154}
\jmlryear{2021}
\jmlrworkshop{LIDTA 2021}

\begin{document}

\maketitle
\begin{abstract}
The ability to collect and store ever more massive databases has been accompanied by the need to process them efficiently.  In many cases, most observations have the same behavior, while a probable small proportion of these observations are abnormal. Detecting the latter, defined as outliers, is one of the major challenges for machine learning applications ($\eg$ in fraud detection or in predictive maintenance). In this paper, we propose a methodology addressing the problem of outlier detection, by learning a data-driven scoring function defined on the feature space which reflects the degree of abnormality of the observations. This scoring function is learnt through a well-designed binary classification problem whose empirical criterion takes the form of a two-sample linear rank statistics on which theoretical results are available. We illustrate our methodology with preliminary encouraging numerical experiments.

\end{abstract}


\section{Introduction}
The problem of ranking multivariate data by degree of abnormality,  referred to as \textit{anomaly ranking}, is of central importance for a wide variety of applications ($\eg$ fraud detection, fleet monitoring, predictive maintenance).  In the standard setup, the 'normal' behavior of the system under study (in the sense of 'not abnormal', without any link to the Gaussian distribution) is described by the (unknown) distribution $F(dx)$ of a generic $\rv$ $X$, valued in $\mathbb{R}^d$. The goal pursued is to build a scoring function $s:\mathbb{R}^d\rightarrow \mathbb{R}_+\cup \{+\infty\}$ that ranks any observations $x_1,\; \ldots,\; x_n$ nearly in the same order as any increasing transform of the density $f$ would do. Ideally, the smaller the score $s(x)$ of an observation $x$ in $\mathbb{R}^d$, the more abnormal it should be considered. In \cite{CLemThom}, a functional criterion, namely a Probability-Measure plot referred to as the \textit{Mass-Volume} curve (the $\mv$ curve in abbreviated form), has been proposed to evaluate the anomaly ranking performance of any scoring rule $s(x)$. This performance measure can be viewed as the unsupervised version of the \textit{Receiver Operating Characteristic} ($\roc$) curve, the gold standard measure to evaluate the accuracy of scoring functions in the bipartite ranking context, see \textit{e.g.} \cite{CV09ieee}. Beyond this approach, let us highlight that the problem of anomaly detection has also been studied \textit{via} various other modelings. For instance, the works of \cite{BergmanHosh20} and \cite{steinwart05a} are based on classification methods, while \cite{Liu2008} build on peeling, \cite{Breunig2000LOF} on local averaging criteria, \cite{FreHab17} on ranking and \cite{Scholkopf2001} on plug-in techniques.

In this paper, we propose a novel two-stage method for detecting and ranking abnormal instances, by means of scalar criteria summarizing the $\mv$ curve and extending the area under its curve, when $F(dx)$ has compact support. Briefly, starting from a sample of observations $X_1,\; \ldots,\; X_n$, we artificially generate an independent second sample $U_1, \ldots, U_m$ that is used as a proxy for outliers. For theoretical reasons explained in the paper, the agnostic choice consists in sampling the $U_i$'s  $\iid$ from the uniform law on a subset of $\mathbb{R}^d$, which $F(dx)$'s support is supposedly included in. 
We then learn to discriminate the $X_i$'s from the $U_i$'s thanks to a scoring function that maximizes two-sample empirical counterparts of the aforementioned criteria, that are in particular robust to imbalanced datasets. The resulting scoring function allows to rank the $X_i$'s by degree of abnormality. 
This novel class of criteria is based on theoretical guarantees provided by \cite{CleLimVay21} on general classes of two-sample linear rank processes, that  incidentally circumvent the difficulty of optimizing the functional $\mv$ criterion. Beyond the classical results of statistical learning theory for these processes, \cite{CleLimVay21} obtain theoretical generalization guarantees for their empirical optimizers. The numerical results performed at the end of the paper also provide strong empirical evidence of the relevance of the approach promoted here.  
\medskip

The article is structured as follows. In section \ref{sec:mot_prelim}, the formulation of the (unsupervised) anomaly ranking problem is recalled at length, together with the concept of $\mv$ curve. In section \ref{sec:method}, the anomaly ranking performance criteria proposed are introduced and their statistical estimation is discussed. Optimization of the statistical counterparts of the criteria introduced to build accurate anomaly scoring functions is also put forward therein. Finally, the relevance of this approach is illustrated by numerical results in section \ref{sec:expes}.

\section{Background and Preliminaries}\label{sec:mot_prelim}

We start off with recalling the formulation of the (unsupervised) anomaly ranking problem and introducing notations that shall be used here and throughout. By $\lambda$ is meant the Lebesgue measure on $\RR^d$, by $\mathbb{I}\{ \mathcal{E} \}$ the indicator function of any event $\mathcal{E}$, while the generalized inverse of any cumulative distribution function $K(t)$ on $\RR$ is denoted by $K^{-1}(u)=\inf\{ t\in \RR:\; K(t)\geq u \}$.
We consider a $\rv$ $X$ valued in $\mathbb{R}^d$, $d\geq 1$, with distribution $F(dx)=f(x)\lambda(dx)$, modeling the 'normal' behavior of the system under study. The observations at disposal $X_1,\; \ldots,\; X_n$, with $n\geq 1$, are independent copies of $X$. Based on the $X_i$'s our goal is to learn a ranking rule for deciding among two observations $x$ and $x'$ in $\mathbb{R}^d$ which one is more 'abnormal'.  The simplest way of defining a preorder\footnote{A preorder $\preccurlyeq$ on a set $\Z$ is a reflexive and transitive binary relation on $\Z$. It is said to be \textit{total}, when either $z\preccurlyeq z'$ or else $z'\preccurlyeq z$ holds true, for all $(z,z')\in\Z^2$.} on $\mathbb{R}^d$ consists in transporting the natural order on $\mathbb{R}_+\cup\{+\infty\}$ onto it through a \textit{scoring function}, \textit{i.e.} a Borel measurable mapping $s:\RR^d \rightarrow \RR_+$: given two observations $x$ and $x'$ in $\mathbb{R}^d$, $x$ is said to be more abnormal according to $s$ than $x'$ when $s(x)\leq s(x')$. The set of all anomaly scoring functions that are integrable with respect to Lebesgue measure is denoted by $\mathcal{S}$. The integrability condition is not restrictive since the preorder induced by any scoring function is invariant under strictly increasing transformation (\textit{i.e.} the scoring function $s$ and its transform $T\circ s$ define the same preorder on $\mathbb{R}^d$ provided that the Borel measurable transform $T:  \text{Im}(s)\rightarrow \RR_+$ is strictly increasing on the image of the $\rv$ $s(X)$, denoted by $\text{Im}(s)$).
One wishes to build, from the 'normal' observations only, a scoring function $s$ such that, ideally,  the smaller $s(X)$, the more abnormal the observation $X$. The set of optimal scoring rules in $\mathcal{S}$ should be thus composed of strictly increasing transforms of the density function $f(x)$ that are integrable $\wrt$ to $\lambda$, namely:
\begin{equation}\label{setoptim}
    \mathcal{S}^*=\{T\circ f:\; T:\text{Im}(f)\rightarrow \RR_+\text{ strictly increasing},\; \int_{\mathbb{R}^d}T\circ f(x)\lambda(dx)<+\infty  \}~.
\end{equation}

The technical assumptions listed below are required to define a criterion, whose optimal elements coincide with $\mathcal{S}^*$.

\begin{itemize}
\item[] $\mathbf{H_1}$ The $\rv$ $f(X)$ is continuous, \textit{i.e.} $\forall c\in \RR_+$, $\mathbb{P}\{f(X)=c  \}=0$.
\item[] $\mathbf{H_2}$ The density function $f(x)$ is bounded: $\vert\vert f\vert\vert_{\infty}\overset{def}{=}\sup_{x\in \mathbb{R}^d}\vert f(x)\vert<+\infty$.
\end{itemize}

\medskip

\begin{figure}[ht]
 \begin{center}
\includegraphics[scale=0.55]{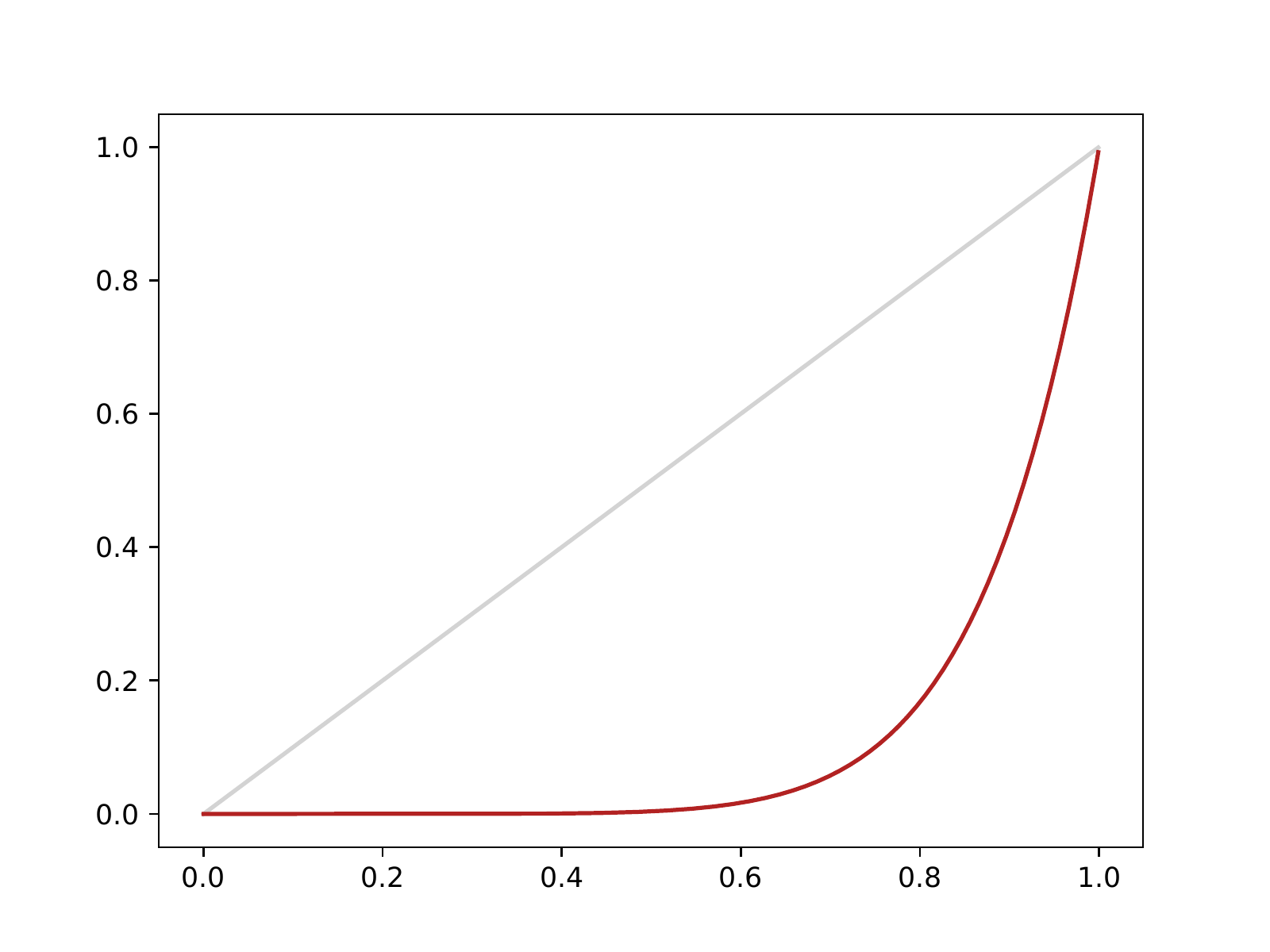}
 \caption{Typical $\mv$ curve in red ($x$-axis:volume, $y$-axis:mass). In gray, the diagonal $y=x$.}
 \label{fig:MVcurve}
 \end{center}
 \end{figure}

 \paragraph{Measuring anomaly scoring accuracy - The $\mv$ curve.}
Consider an arbitrary scoring function $s\in \mathcal{S}$ and denoted by $\Omega_{s,t}=\{x\in \X:\; s(x)\geq t \}$, $t\geq 0$, its level sets. As $s$ is $\lambda$-integrable, the measure $\lambda(\Omega_{s,t})\leq (\int_{u\in \RR_+}s(u)du)/t$ is finite for any $t>0$. Introduced in \cite{CLemThom}, a natural measure of the anomaly ranking performance of any scoring function candidate $s$ is the Probability-Measure plot, referred to as the \textit{Mass-Volume} ($\mv$) curve:
\begin{equation} \label{MVtheo1}
t>0 \mapsto \left(\mathbb{P}\{s(X)\geq t\},\; \lambda(\{x\in \mathbb{R}^d:\; s(x)\geq t \})  \right)=\left(F(\Omega_{s,t}),\; \lambda(\Omega_{s,t})  \right)~.
\end{equation}
 Connecting points corresponding to possible jumps, this parametric curve can be viewed as the plot of the continuous mapping $\mv_s:\alpha\in (0,1)\mapsto \mv_s(\alpha)$, starting at $(0,0)$ and reaching $(1,\; \lambda\bigl(\supp(F)\bigr)$ in the case where the support $\supp(F)$ of the distribution $F(dx)$ is compact, or having the vertical line '$\alpha=1$' as an asymptote otherwise. A typical $\mv$ curve is depicted in Fig. \ref{fig:MVcurve}.

 Let $\alpha\in (0,1)$. Denoting by $F_s(t)$ the cumulative distribution function of the $\rv$ $s(X)$, we have:
\begin{equation} \label{MVtheo2}
\text{MV}_s(\alpha)=\lambda\left(\{ x \in \mathbb{R}^d:\; s(x)\geq F_s^{-1}(1-\alpha) \}  \right),
\end{equation}
when $F_s\circ F_s^{-1}(\alpha)=\alpha$.  This functional criterion is invariant by increasing transform and induces a partial order over the set $\S$. Let $(s_1,s_2)\in \mathcal{S}^2$, the ordering defined by $s_1$ is said to be more accurate than the one induced by $s_2$ when:
 \begin{equation*}
     \forall \alpha \in (0,1),\;\; \mv_{s_1}(\alpha)\leq \mv_{s_2}(\alpha)~.
 \end{equation*}
 
 As summarized by the result stated below, the $\mv$ curve criterion is adequate to measure the accuracy of scoring functions with respect to anomaly ranking. 

It reveals in particular that optimal scoring functions are those whose $\mv$ curve is minimum everywhere.

 \begin{proposition} \label{prop:opt} (\cite{CLemThom}) Let the assumptions $\mathbf{H_1}-\mathbf{H_2}$ be fulfilled.
The elements of the class $\mathcal{S}^*$ have the same (convex) $\mv$ curve and provide the best possible preorder on $\mathbb{R}^d$ w.r.t. the $\mv$ curve criterion:
\begin{equation}\label{eq:opt}
\forall (s,\alpha)\in \mathcal{S}\times (0,1),\;\; \mv^*(\alpha)\leq \mv_s(\alpha)~,
\end{equation}
where $\mv^*(\alpha)=\mv_{f}(\alpha)$ for all $\alpha\in (0,1)$.

 \end{proposition}

Equation \eqref{eq:opt} reveals that the lowest the $\mv$ curve (everywhere) of a scoring function $s(x)$, the closer the preorder defined by $s(x)$ is to that induced by $f(x)$.
Favorable situations are those where the $\mv$ curve increases slowly and rises more rapidly when coming closer to the 'one' value: this correponds to the case where $F(dx)$ is much concentrated around its modes, $s(X)$ takes its highest values near the latter and its lowest values are located in the tail region of the distribution $F(dx)$. Incidentally, observe that the optimal curve $\mv^*$ somehow measures the spread of the distribution $F(dx)$ in particular for 
large values of $\alpha$ $\wrt$ extremal observations (\textit{e.g.} a light tail behavior corresponds to the situation where $\mv^*(\alpha)$ increases rapidly when approaching $1$), whereas it should be examined for small values of $\alpha$ when modes of the underlying distributions are investigated (a flat curve near $0$ indicates a high degree of concentration of $F(dx)$ near its modes).
\medskip

{\bf Statistical estimation.} In practice, the $\mv$ curve of a scoring function $s\in \mathcal{S}$ is generally unknown, just like the distribution $F(dx)$, and it must be estimated. A natural empirical counterpart can be obtained by plotting the stepwise graph of the mapping:
\begin{equation}\label{eq:emp_curve}
\widehat{\mv}_s(\alpha) :\alpha\in (0,1)\mapsto \lambda\left(\left\{x\in \mathbb{R}^d:\; s(x)\geq \widehat{F}_{s,n}^{-1}(1-\alpha)  \right\}  \right)~,
\end{equation}
where $\widehat{F}_{s,n}(t)=(1/n)\sum_{i=1}^n\mathbb{I}\{s(X_i)\leq t  \}$ denotes the empirical $\cdf$ of the $\rv$ $s(X)$ and $\widehat{F}_{s,n}^{-1}$ its generalized inverse. In \cite{CLemThom}, for a fixed $s\in \mathcal{S}$, consistency and asymptotic Gaussianity (in $\sup$-norm)
of the estimator \eqref{eq:emp_curve} has been established, together with the asymptotic validity of a smoothed bootstrap procedure to build confidence regions in the $\mv$ space. However, depending on the geometry of the superlevel sets of $s(x)$, it can be far from simple to compute the volumes. In the case where $F$ has compact support, included in $[0,1]^d$ say for simplicity, and from now on it is assumed it is the case, they can be estimated by means of Monte-Carlo simulation. Indeed, if one generates a synthetic $\iid$ sample $\{U_1, \ldots, U_m\}$, independent from the $X_i$'s and drawn from the uniform distribution on $[0,1]^d$, which we denote by $\mathcal{U}_{d}$, a natural estimator of the volume $\widehat{\mv}_s(\alpha)$ is:
\begin{equation}
    \widetilde{\mv}_s(\alpha)=\frac{1}{m}\sum_{j=1}^m\mathbb{I}\{s(U_j)\geq \widehat{F}_{s,n}^{-1}(1-\alpha)  \}~.
\end{equation}

\medskip

\paragraph{Minimization of the empirical area under the $\mv$ curve.} Thanks to the $\mv$ curve criterion, it is possible to develop a statistical theory for the anomaly scoring problem. From a statistical learning angle, the goal is to build from training data $X_1,\; \dots,\; X_n$ a scoring function with $\mv$ curve as close as possible to $\mv^*$. Whereas the closeness between (continuous) curves can be measured in many ways, the $L_1$-distance offers crucial advantages. Indeed, we have:
 \begin{equation*}
 d_1(s,f)=\int_{\alpha=0}^{1}\vert \mv_s(\alpha)-\mv^*(\alpha) \vert \, d\alpha=\int_{\alpha=0}^{1} \mv_s(\alpha)d\alpha-\int_{\alpha=0}^{1}\mv^*(\alpha)d\alpha ~,\\
 \end{equation*}
Notice that $d_{1}(s,f)$, $i\in \{1,\; \infty  \}$, is not a distance between the scoring functions $s$ and $f$ but measures the dissimilarity between the preorders they define and that minimizing $d_1(s,f)$ boils down to minimizing the scalar quantity $\int_{\alpha=0}^{1-\varepsilon} \mv_s(\alpha)d\alpha$, the area under the $\mv$ curve. From a practical perspective, one may then learn an anomaly scoring rule by minimizing the empirical quantity:
\begin{equation*}
    \int_0^1 \widetilde{\mv}_s (\alpha)  d \alpha~.
\end{equation*}

This boils down to maximizing the rank-sum (or Wilcoxon Mann-Whithney) statistic (see \cite{Wil45}) given by:
 \begin{equation}\label{eq:ranksum}
  \widehat{W}_{n,m}(s)=\sum_{i=1}^n \Rank(s(X_i))~,
 \end{equation}
 where $\Rank(s(X_i))$ is the rank of $s(X_i)$ among the pooled sample $\{s(X_1),\; \ldots,\; s(X_n)\}\cup\{s(U_1),\; \ldots,\; s(U_m)\}$: $\Rank(s(X_i))= \sum_{l=1}^n\mathbb{I}\{s(X_l)\leq s(X_i)  \} +  \sum_{j=1}^m\mathbb{I}\{s(U_j)\leq s(X_i)  \}$. Indeed, just like the empirical area under the $\roc$ curve  can be related to the rank-sum statistic, we have:
\begin{equation}\label{eq:mvtoMWW}
 nm\left(1-\int_0^1\widetilde{\mv}_s(\alpha)d\alpha\right)+n(n+1)/2=\widehat{W}_{n,m}(s)~.
\end{equation}
In the next section, we introduce more general empirical summaries of the $\mv$ curve that are of the form of two-sample rank statistics, just like \eqref{eq:ranksum}, and propose to solve the anomaly ranking problem through the  maximization of the latter.

\section{Measuring and Optimizing Anomaly Ranking Performance}\label{sec:method}

In this section, a class of anomaly ranking performance criteria are introduced, which can be estimated by two-sample rank statistics. We also emphasize that a natural approach to anomaly ranking consists in maximizing such empirical scalar criteria.

\subsection{Scalar Criteria of Performance and Two-sample Rank Statistics}\label{sub:meth_pbform}

Here we develop the statistical learning framework we propose for anomaly ranking. Let $p\in (0,1)$, we assume that $N\geq 2$ observations are available: $n = \lfloor pN \rfloor$ 'normal' $\iid$ observations $X_1, \ldots, X_n$ taking their values in $[0,1]^d$ for simplicity drawn from $F(dx)=f(x)\lambda(dx)$ and $m = N - n$ $\iid$ realizations of the uniform distribution $\mathcal{U}_{d}$, independent from the $X_i$'s. Hence, $p$ represents the 'theoretical' proportion of 'normal' observations among the pooled sample. Let a class of scoring functions $\S_0 \subset \S$ such that, for all $s(x)$, we consider the mixture distribution $G_s = pF_s + (1-p) \lambda_s$ and its empirical counterpart  $\widehat{G}_{s,N}(t)=(1/n)\sum_{i=1}^n\mathbb{I}\{s(X_i)\leq t  \} + (1/m)\sum_{j=1}^m\mathbb{I}\{s(U_i)\leq t  \}$.
Notice that since $n/N\rightarrow p$ as $N$ tends to infinity, the quantity above is a natural estimator of the \cdf $G_s$. We refer to the \textit{scored} random samples for $\{s(X_1), \ldots, s(X_n)\}$ and $\{s(U_1), \ldots, s(U_m)\}$. Therefore, motivated by Eq. \eqref{eq:mvtoMWW}, Definition \ref{def:Wphi} below provides the class of \textit{ $W_{\phi}$-performance criteria} we consider in the subsequent procedure. 
\begin{definition}\label{def:Wphi} Let $\phi:[0,1]\rightarrow  \RR$ be a nondecreasing function. The '$W_{\phi}$-ranking performance criterion' with 'score-generating function' $\phi(u)$ based on the mixture cdf $G_s(dt)$ is given by:
	\begin{equation}\label{eq:Wtrue}
		W_{\phi}(s) = \mathbb{E} [(\phi \circ G_s)(s(X))]~.
	\end{equation}
\end{definition}

One can naturally relate this generalized form to the $\mv$ curve, justifying this choice of scalar performance criteria as summaries of the $\mv$ curve, through the equality:

\begin{equation}\label{eq:wphitomv}
   W_{\phi}(s) = \int_{0}^1 \phi\left(1 -  p \alpha -(1-p)\mv_s (\alpha) \right) d\alpha~.
\end{equation}

Equipped with the two random samples, the following Definition \ref{def:R_stat} provides an empirical counterpart, that generalizes the empirical summaries of the $\mv$ curve \textit{via} collections of two-sample linear rank statistics. Precisely, for a given mapping $s(x)$, we allow to weight the sequence of 'normal ranks' $\ie$ the ranks of the scored 'normal' instances among the pooled sample, by means of a \textit{score-generating} function.


\begin{definition}{\sc (Two-sample linear rank statistics)}\label{def:R_stat} Let $\phi:[0,1]\rightarrow  \RR$ be a nondecreasing function. The two-sample linear rank statistics with 'score-generating function' $\phi(u)$ based on the random samples $\{X_1,\;\ldots,\; X_n\}$ and  $\{U_1,\;\ldots,\; U_m\}$ is given by:
	\begin{equation}\label{eq:R_stat}
		\widehat{W}^{\phi}_{n,m}(s)=\sum_{i=1}^n\phi\left( \frac{\Rank(s(X_i))}{N+1} \right)~,
	\end{equation}
	where $\Rank(t) = N \widehat{G}_{s,N}(t) =  \sum_{i=1}^n\mathbb{I}\{s(X_i)\leq t  \} +  \sum_{j=1}^m\mathbb{I}\{s(U_j)\leq t  \}$. 
\end{definition}

\paragraph{Optimality.}

Briefly, we refer to the comprehensive analysis of the general class of criteria in \cite{CleLimVay21}, that establishes the theoretical guarantees for the consistency of the two-stage procedure we detail in the following subsection. 
Importantly, the set of optimal maximizers of the empirical $W_{\phi}$-criteria coincides with the nondecreasing transforms of the likelihood ratio, just like for the $\mv$ curves, as shown thourgh the Eq. \eqref{eq:wphitomv}.

The optimal set $\S^*$ derived in Eq. \eqref{setoptim} underlines the implicit characterization that inherits an outlier: the lower the scalar score is and the likelier anomalous the observation can be considered. Also, the notion of distance induced by the rank-based criteria is in fact directly related to the distribution of the 'normal' sample compared to the Uniform one.

\medskip

\paragraph{Choosing $\phi$.} As foreshadowed above, the choice of the score-generating function is an asset of this class of criteria as it provides a flexibility $\wrt$ the weighting of the area under the $\mv$ curve. Indeed, its minimization directly implies the maximization of the $W_{\phi}$-criterion (see Eq. \eqref{eq:wphitomv}), recalling the nondecreasing variation of $\phi(u)$. Therefore, one can hope to recover at best the $\mv^*$ curve by the right choice of $\phi(u)$, especially when the initial sample is noisy. 
Additionally, when going back to the problem of learning to rank the (possible abnormal) instances, it is an advantage to weight the ranks accordingly. 

\begin{figure}[ht!]
    \centering
    \includegraphics[scale=0.3]{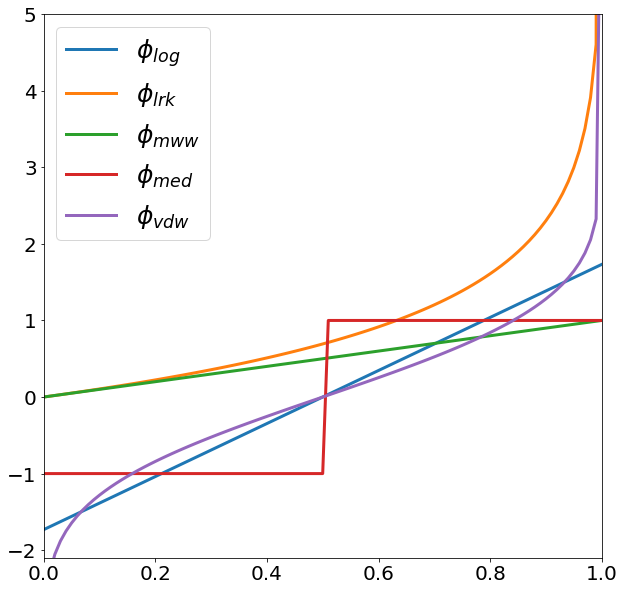}
    \caption{Curves of two-sample score-generating functions with the associated statistical test: Logistic test $\phi_{log}(u)= 2 \sqrt{3}(u - 1/2)$ in blue,  Logrank test $\phi_{lrk}(u)= - \log(1-x)$ in orange, Mann-Whitney-Wilcoxon test $\phi_{mww}(u)= u$ in green,  Median test  $\phi_{med}(u)= \text{sgn} (u - 1/2)$  in red,  Van der Waerden test $\phi_{vdw}(u)= \Phi^{-1}(u)$  in purple, $\Phi$ being the normal quantile function.}
    \label{fig:scoregenstatclass}
\end{figure}

First, we recall the simplest uniform weighting of each 'normal' rank with $\phi(u)=u$. It parenthetically yields to Eq. \eqref{eq:mvtoMWW}, of continuous version:  $W(s) =p/2 + (1-p)(1-\int_{0}^1 \mv_s (\alpha) d\alpha)$, where the area under the $\mv$ curve is clearly computed. 
  Other functions were introduced in the literature related to classic univariate two-sample rank statistics. Figure \ref{fig:scoregenstatclass} gathers classical nondecreasing score-generating functions broadly used for two-sample statistical tests (refer to \cite{Haj62}).

\subsection{The Two-Stage Procedure}\label{sub:meth_proc}

In this paragraph, we detail the two-stage procedure, where we assume that both the framework and assumptions detailed in the previous subsection are adopted. We define the test sample as the set of $\iid$ random variables $\{X_1^t,\;\ldots,\; X_{n_t}^t\}$, with $n_t\in \NN^*$, \textit{a priori} drawn from $F(dx)$. The goal pursued is to distinguish among the test sample, the instances the most likelier to be anomalous. In particular, we propose a first step $(1.)$ that outputs an optimal ranking rule $\hat{s}_{n,m}(x)$, in the sense of the maximization of the rank statistics of Eq. \eqref{def:R_stat}. Then, in the second step $(2.)$ and equipped with this rule, the instances of the test sample are optimally ranked by increasing order of similarity \wrt the $X$'s. We also choose to watch a number of $n_{lowest} \in \NN^*$ worst ranked instances \ie of lowest empirical score. The procedure is detailed in the following Fig. \ref{fig:rankprocedure}. By means of the recalled theoretical guarantees proved in \cite{CleLimVay21}, it results to the asymptotic consistency of step $(1.)$ as well as its nonasymptotic consistency with high probability, under some technical assumptions.

\begin{figure}[ht!]
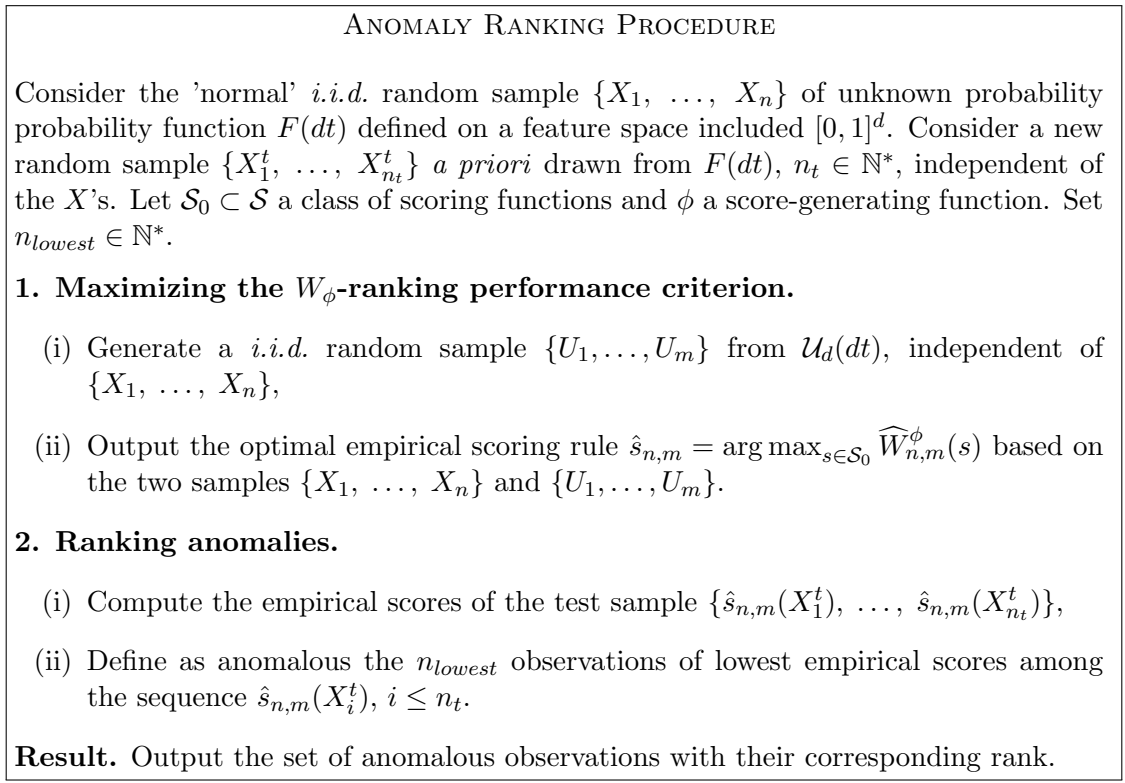

    \centering
\fbox{\begin{minipage}{0.95\textwidth} \label{twostageproc}
		
		\begin{center}{\sc Anomaly Ranking Procedure}\end{center}

		Consider the 'normal' $\iid$ random sample $\{X_1,\;\ldots,\; X_n\}$ of unknown probability probability function $F(dt)$ defined on a feature space included $[0,1]^d$. Consider a new random sample $\{X_1^t,\;\ldots,\; X_{n_t}^t\}$  \textit{a priori} drawn from $F(dt)$, $n_t\in \NN^*$, independent of the $X$'s. Let $\S_0 \subset \S$ a class of scoring functions and $\phi$ a score-generating function. Set $n_{lowest} \in \NN^*$.
		\medskip
		
		{\bf 1. Maximizing the $W_{\phi}$-ranking performance criterion.}
		\begin{enumerate}
			\item[(i)] Generate a $\iid$ random sample $\{U_1, \ldots, U_m\}$ from $\mathcal{U}_{d}(dt)$, independent of $\{X_1,\;\ldots,\; X_n\}$,
			\item[(ii)] Output the optimal empirical scoring rule $\hat{s}_{n,m} = \argmax_{s \in \S_0} \widehat{W}^{\phi}_{n,m}(s)$ based on the two samples $\{X_1,\;\ldots,\; X_n\}$ and  $\{U_1, \ldots, U_m\}$.
		\end{enumerate}
	
		{\bf 2. Ranking anomalies.} 
		\begin{enumerate}
		\item[(i)]  Compute the empirical scores of the test sample $\{\hat{s}_{n,m} (X_1^t),\;\ldots,\; \hat{s}_{n,m} (X_{n_t}^t)\}$,
		\item[(ii)] Define as anomalous the $n_{lowest}$ observations of lowest empirical scores among the sequence $\hat{s}_{n,m} (X_i^t) $, $i \leq n_t$.
		\end{enumerate}
		{\bf Result.} Output the set of anomalous observations with their corresponding rank.

\end{minipage} }

    \caption{Two-stage procedure for learning to rank anomalies.}
    \label{fig:rankprocedure}
\end{figure}

\section{Numerical Experiments}\label{sec:expes}

In this section, we illustrate the procedure promoted along the paper through numerical experiments on imbalanced synthetic data. As these experiments are mainly here to support our methodology, we propose for the step $(1.)$ to learn the empirical maximizer $\hat{s}_{n,m}$ by means of a regularized classification algorithm. At a technical level, we would ideally like to replace usual loss criterion such as the BCE (Binary Cross-Entropy) loss by our tailored objective $W_\phi$. Unfortunately, the latter is not smooth and of highly correlated terms, which results in many challenges regarding its optimization. In order to incorporate $W_\phi$ and still keeping good performances, we {\it (i)} use a regularized proxy of it and {\it (ii)} incorporate the regularized criterion in a penalization term. The second point allows to drive the learning with a usual BCE loss, which asymptotically amounts to estimate the conditional probability $\mathbb{P}(y=1 \, | \, X)$, while considering $W_\phi$.

\medskip

\begin{figure}[ht!]
	\centering
	\begin{tabular}{cc}
		
		\parbox{7cm}{	
         \includegraphics[scale=0.4]{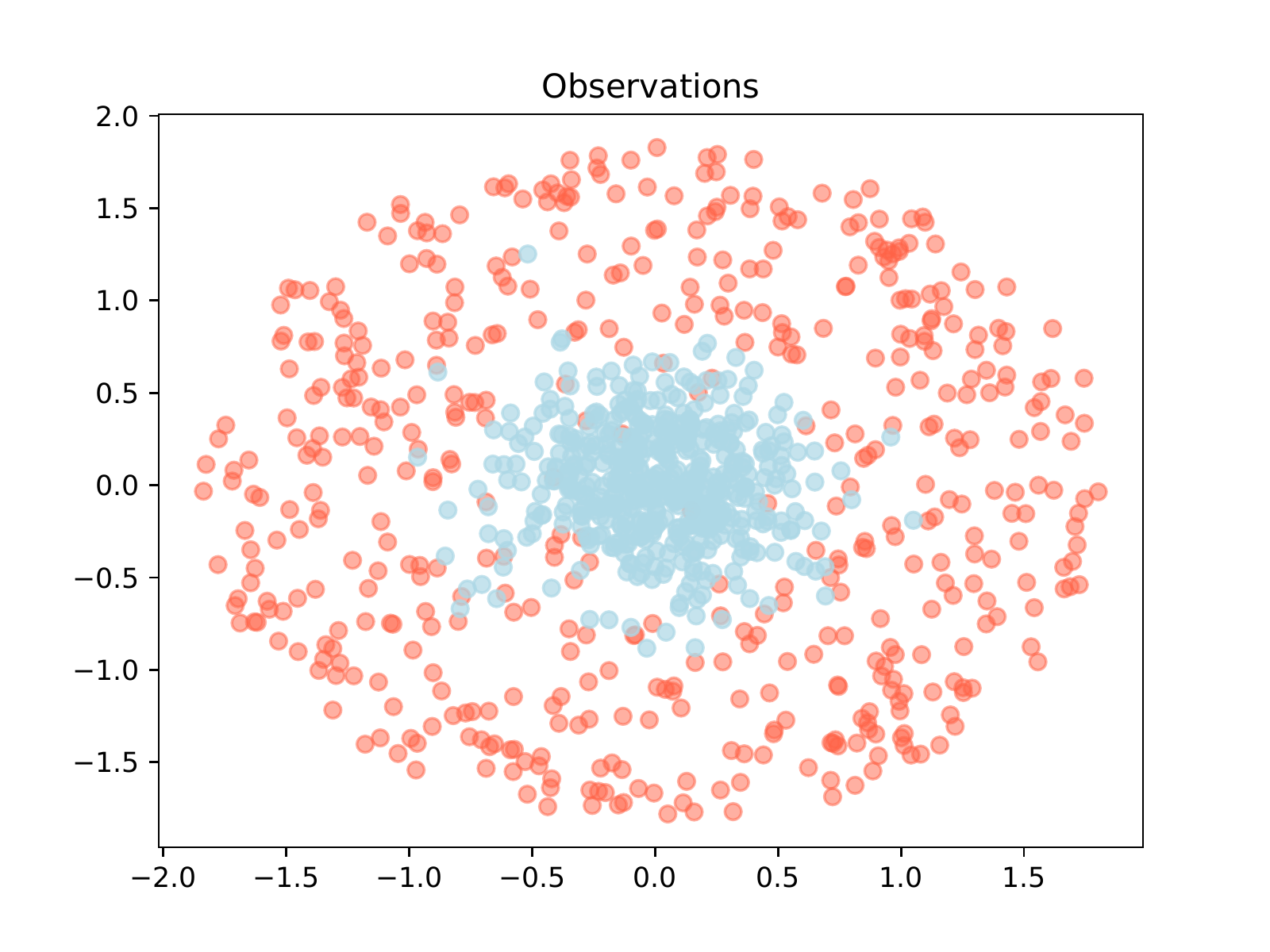}\\
         \small{(a) Train data. $(n, \; m) = (1000, \; 500)$.}
     }
   	\parbox{7cm}{	
         \includegraphics[scale=0.4]{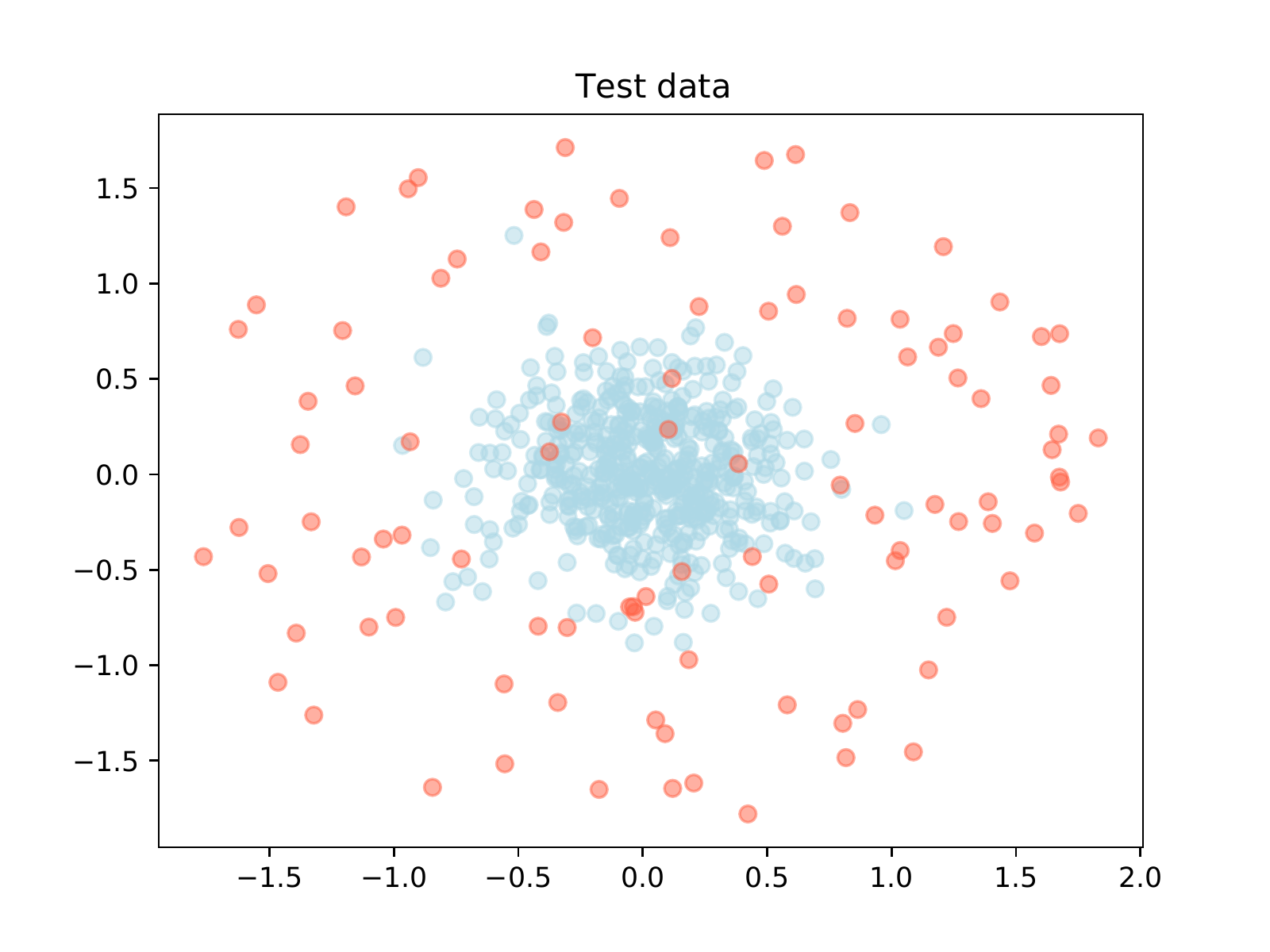}\\
         \small{(b) Test data. $(n_t, \; m_t) = (400, \; 100)$.}
  }
\end{tabular}
        \caption{Data visualization for the two generating processes.  The Gaussian observations are represented in blue. The 'synthetic outliers' samples drawn from the radial law are represented in red. The left figure $(a)$ corresponds to the train dataset, the right $(b)$ to the test dataset. }
        \label{fig:data_gen}
\end{figure}

\paragraph{Data generating process.} We generated the 'positive' sample by $\iid$ Gaussian variables $X_1, \ldots, X_{n}$, $n =1000$, in dimension $d=2$, centered and with covariance matrix $0.1\times I_2$ (where $I_2$ is the identity matrix). We chose the Gaussian law for its attractive structure and in particular for its symmetry, it can be a reasonable choice in many situations where the data at hand are indeed well structured. 
We then sampled the 'negative' sequence of $\iid$ $\rv$ $U_1', \ldots, U_{m}'$, $m=500$, from the following radial law, expressed in terms of its density in polar coordinates: 
\begin{equation*}
     \mathrm{RadLaw}_{\alpha, \beta}: (v,r) \in \mathbb{S}^{d-1} \times (0,1)\mapsto \frac{1}{\mathrm{Area}( \mathbb{S}^{d-1})} dv \times \frac{1}{B(\alpha, \beta)}r^{\alpha-1} (1-r)^{\beta-1} dr~,  
\end{equation*}

where $\alpha, \beta > 0$ are two tunable parameters, $\mathbb{S}^{d-1} = \{ x \in \RR^d, \; \| x \| = 1 \}$ is the unit sphere, and where $B(\alpha, \beta) = \int_0^1 r^{\alpha-1} (1-r)^{\beta-1} dr $. In other words, $v$ is uniformly sampled in the unit sphere and $r$ has Beta law with parameters $\alpha$ and $\beta$. Notice that $\alpha = \beta = 1$ corresponds to the Uniform law and that, when $\beta=1$, the law puts more mass around $1$ as $\alpha>1$ increases. In our experiment, we choose $\alpha=3$ and $\beta =1$. Denoting by $\mathrm{rad} = \max_{1 \leq i \leq n} || X_i ||$, we finally obtained $m$ 'synthetic outliers' $U_1, \ldots, U_{m}$ defined by $U_i = (\mathrm{rad} + \varepsilon) \times U_i'$, with $\varepsilon = 0.01$. To simplify the notations, we denote by $\mathbf{Z}_{train}$ the concatenation of the $X_i$'s and the $U_i$'s. We also denote by $\mathbf{y}_{train}$ the labels, where we choose to assign the label $1$ ($\resp$ $0$) to the 'positive' ($\resp$ 'negative') sample. Figure \ref{fig:data_gen} illustrates both data generating processes. For the test set, we generated similarly a sequence of $n_{t}=400$ $\iid$ Gaussian $\rv$ $X_1^t, \ldots, X_{n_{t}}^t$ from the same Gaussian law as the 'positive' sample, and a $\iid$ random sequence $U_1^{t}, \ldots, U_{m_{t}}^t$, $m_{t}=100$, drawn from the law $\mathrm{RadLaw}_{\alpha_t, \beta_t}$, with $\alpha_t = 2.$ and $\beta_t = 1.$, dilated by a factor $(\mathrm{rad} + \varepsilon)$. 

\medskip

\paragraph{Metrics.} Once the algorithm that learns a (renormalized) optimal scoring function $\hat{s}_{n,m}: \mathbb{R}^d \rightarrow (0,1)$ has been trained ($\ie$ step $(1.)$), we score the test data with $\hat{s}_{n,m}$ and compute the proportion of true outliers among the $n_{lowest}$ points having lowest scores ($\ie$ step $(2.)$). We let $n_{lowest}$ varies in $\{25, 50, 75, 100\}$. Formally, if $\xi_1 \preccurlyeq \cdots \preccurlyeq \xi_{n_{t} +m_{t}}$ denote the points $X_i^t$ and $U_i^t$ sorted by scores, $\ie$ the ordered sequence based on $\hat{s}_{n,m}(\bZ_{test, 1}), \cdots,  $ $\hat{s}_{n,m}(\bZ_{test, n_t + m_t}) $, we compute the following accuracy:
\begin{equation}
 \mathrm{Acc}_{n_{lowest}} = \frac{1}{n_{lowest}} \sum\limits_{i=1}^{n_{lowest}} \mathbb{I}\{ \xi_i \in \{ U_1^t, \ldots, U_{m_{t}}^t \} \}~.
\end{equation}

\medskip 

\paragraph{Neural Network.} We trained a neural network $\textsc{mlp}$ composed of one hidden layer of size $2\times d$, a ReLu activation function and whose last layer is a Sigmoid function, computing the desired score. For each $n_{epoch} = 30$ epochs, we use the following training scheme:
\begin{enumerate}
    \item Each sample of $(\mathbf{Z}_{train},\;  \mathbf{y}_{train})$ is individually passed through the network, the BCE loss is computed\footnote{Remember it is given by $-y\ln \hat{y} - (1-y) \ln(1- \hat{y})$, where $\hat{y}=\textsc{mlp}(X)$.} and a backpropagation step is performed,
    \item At the end of each epoch, the whole batch of the training dataset $(\mathbf{Z}_{train}, \; \mathbf{y}_{train})$ is passed through the network and we computed the Binary Cross Entropy loss, denoted by $\mathrm{BCE}$, and the following proxy of $W_\phi$:
    \begin{equation*}
        \widehat{W}_{n,m}^{\phi} =  \sum\limits_{i=1}^{n} \phi \left( \frac{ (n+m) \times  \textsc{mlp}(X_i)  + 1 }{n +m +1} \right)  ~.
    \end{equation*}
    In our experiments, we choose $\phi(u)=u$ 
    and $\phi_{u_0}(u)=  u \mathbb{I}\{u\geq u_0\}$ with $u_0=0.7$, as defined in section \ref{sub:meth_pbform}. 
    We then compute the regularized loss $
        \mathrm{BCE} - \lambda \widehat{W}_{n,m}^{\phi}$,
    where $\lambda$ is a hyperparameter in $\{0, 0.01, 0.1, 1, 10\}$.
\end{enumerate}

The training procedure of the Neural Net is summarized in the Algorithm \ref{algo:NN}.

\begin{algorithm2e}[ht!]\label{algo:GA}
	\SetAlgoLined
	
	\KwData{$(\bZ_{train}, \by_{train})$.}
	\KwIn{Network $\textsc{mlp}$, number of epochs $n_{epoch}$, penalization strength $\lambda$.}

	\KwResult{Trained network.}
	\BlankLine
	
	\For{$n=0, \ldots, \; n_{epoch}$}{
		 \For{ $X, y \in \mathbf{Z}_{train}, \mathbf{y}_{train}$}{
		 compute $\hat{y} = \textsc{mlp}$ ; \\
		 compute $BCE = BCE(\hat{y}, y)$, backpropagate and zero\_grad ;
		 } 	
		 compute $\mathbf{\hat{y}} = \textsc{mlp}( \mathbf{Z}_{train})$ ; \\
		 compute $\mathrm{BCE} = \mathrm{BCE}(\mathbf{\hat{y}}, \mathbf{y})$ and $\widehat{W}_{n,m}^{\phi}$ ; \\
		 compute the regularized loss $\mathrm{BCE} - \lambda \widehat{W}_{n,m}^{\phi}$, backpropagate and zero\_grad ;
	}

\caption{Training of the Neural Network}
\label{algo:NN}
\end{algorithm2e}

\smallskip

\paragraph{Repetitions.} We repeat $B = 100$ times the procedure, each time computing the accuracy metric defined above.

\medskip

\paragraph{Visualization and results.} In this section, we only display the results obtained with $\phi(u)=u$ since they are very similar to the one obtained with $\phi(u) = u \mathbb{I}\{ u\geq u_0\}$. This is probably due to the very simple framework adopted for the data generating process and further investigations would be of interest.

\begin{figure}[ht!]
    \hspace{-1.9cm}
    \includegraphics[scale=0.4]{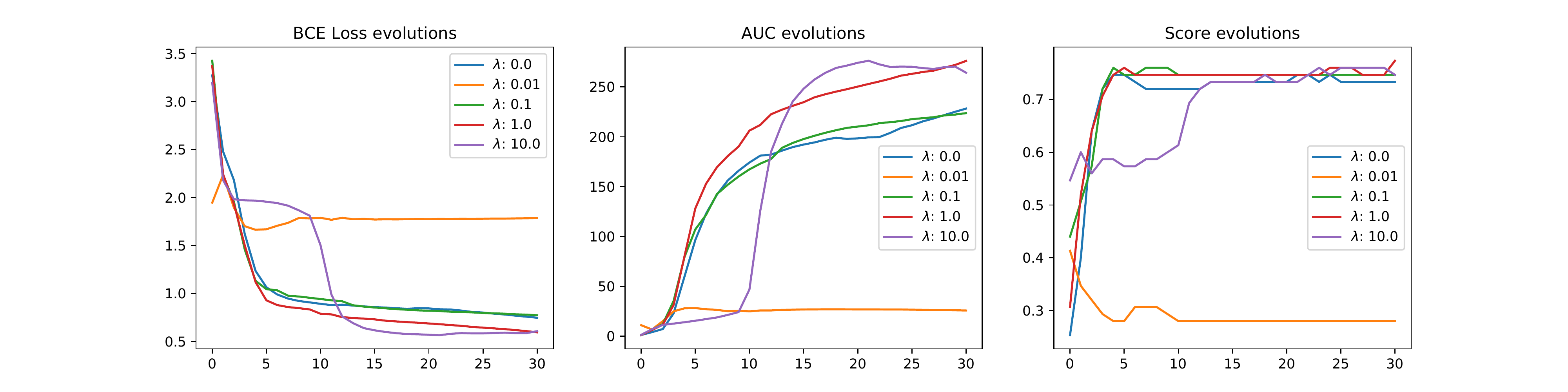}
    \caption{Evolutions of the BCE loss, the AUC proxy and the accuracy for $n_{lowest}=75$ in function of the epochs, for $\phi(u)=u$ and all values of the hyperparameter $\lambda \in \{0, 0.01, 0.1, 1, 10\}$.}
    \label{fig:training}
\end{figure}

For the first learning loop, we saved the evolution of the BCE losses, for all values of $\lambda$, computed at each epoch together with the $W_\phi$ proxy and the accuracy metric for $n_{lowest} = 75$. As displayed in Figure \ref{fig:training}, one can see that the incorporation of the empirical $W_\phi$ criterion in the penalization term improves the performances for a well chosen parameter $\lambda$. For instance, $\lambda \in \{ 1, 10\}$ output the best results in this setting. 

\begin{figure}[ht!]
    \centering
    \includegraphics[scale=0.6]{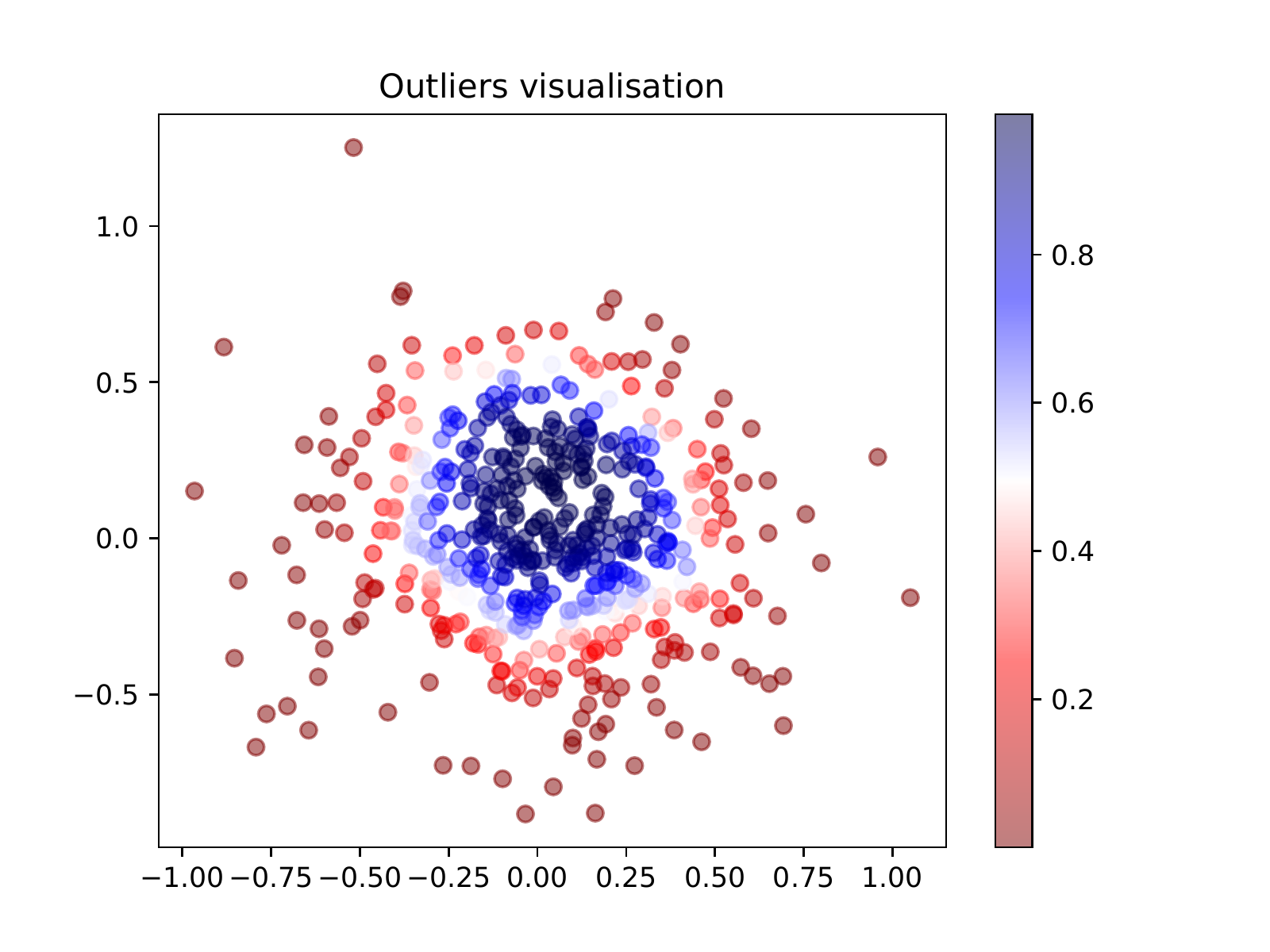}
    \caption{A heatmap of the scores for $\phi(u) = u$.}
    \label{fig:heat_map}
\end{figure}

\begin{figure}[ht!]
	\centering
	\begin{tabular}{cc}
		\parbox{7cm}{	
			\includegraphics[scale=0.4]{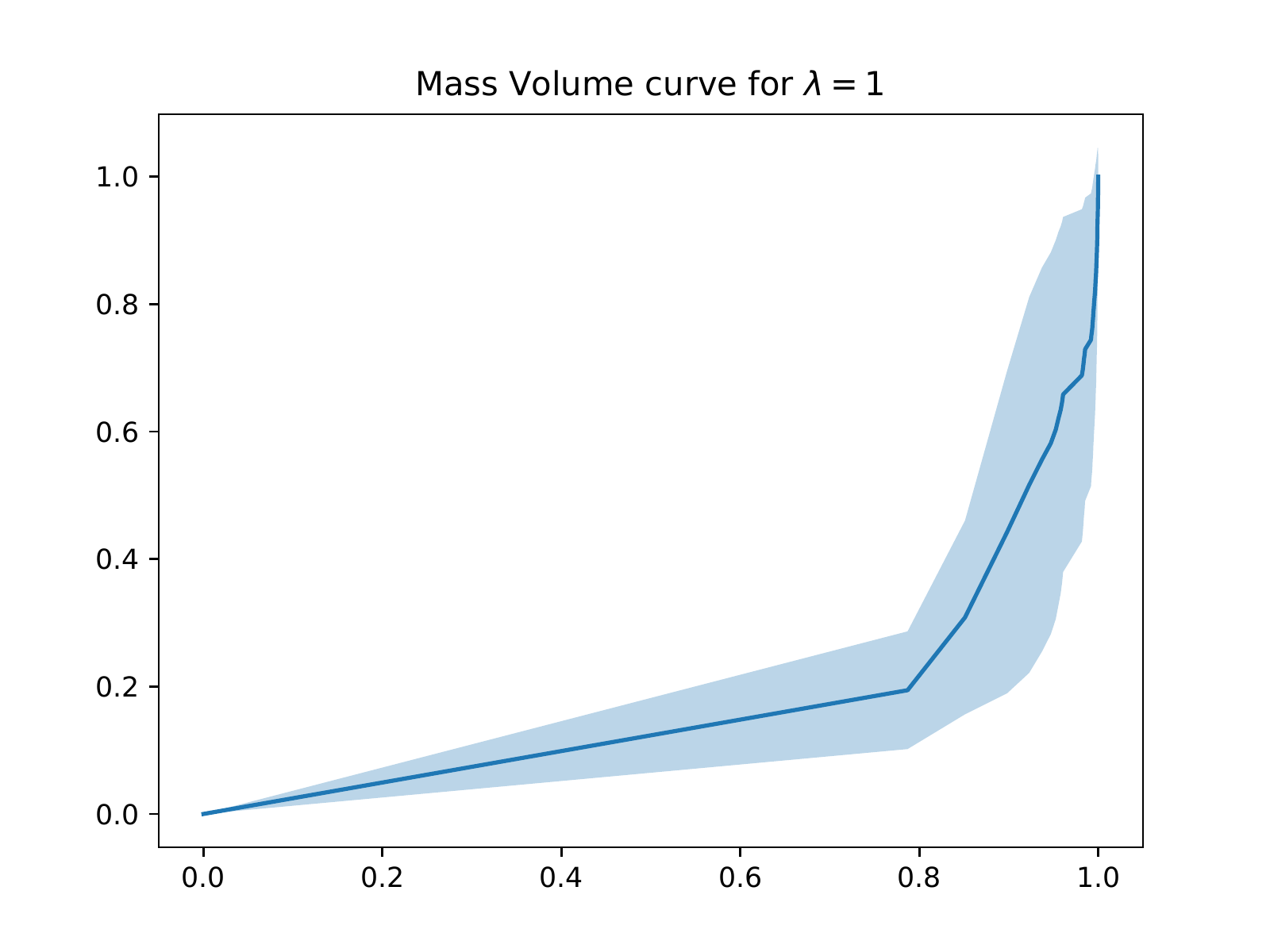}\\
			 \small{(a) $\lambda = 1$ and $\phi(u)=u$.}
		}
		\parbox{7cm}{
			\includegraphics[scale=0.4]{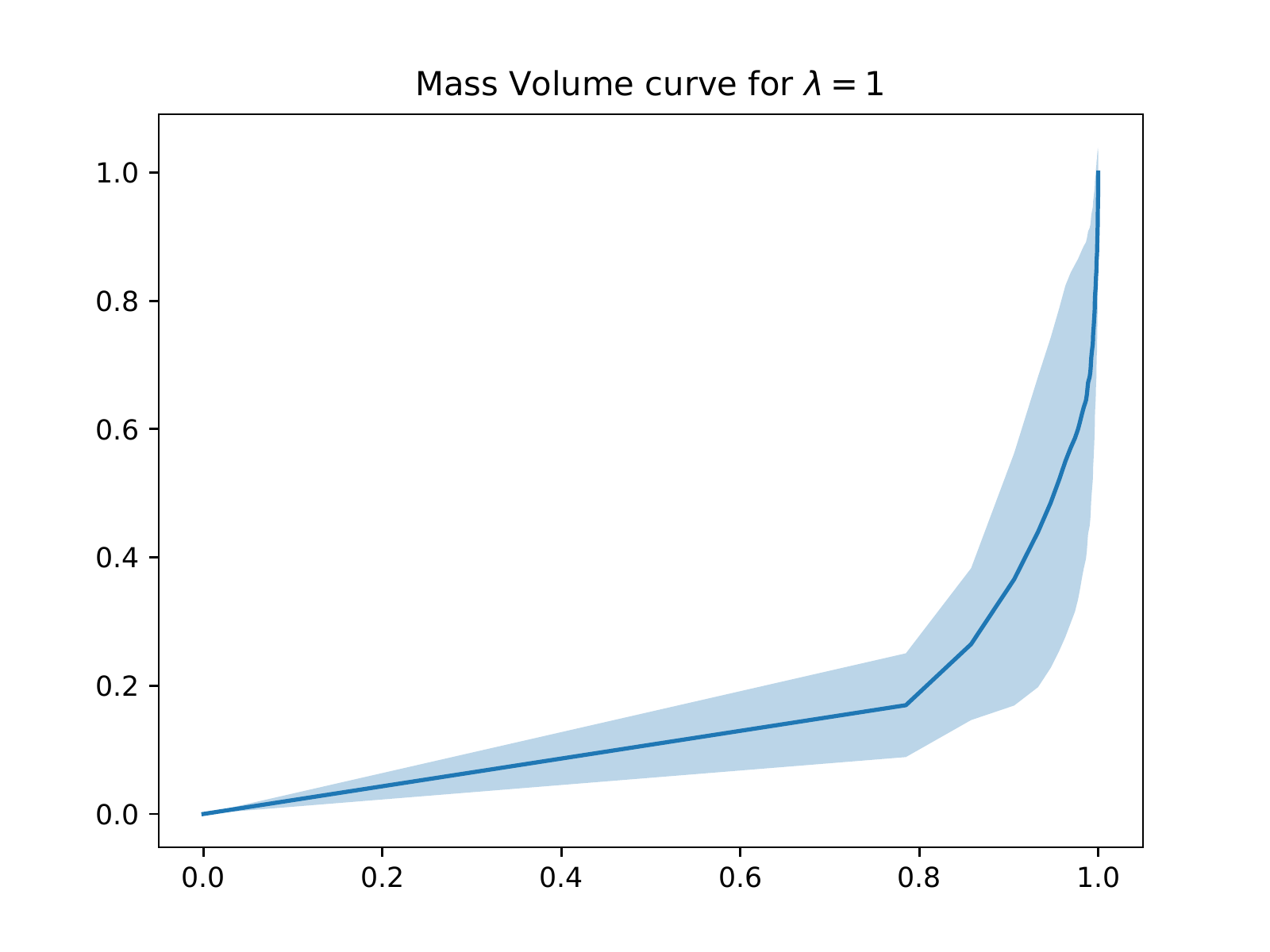}\\
			\small{(b) $\lambda = 1$ and $\phi(u)=u \mathbb{I}\{ u \geq u_0\}$.}
		}
		\medskip
	\end{tabular}
        \caption{Empirical Mass-Volume curves.}
\label{fig:mass_volume}
\end{figure}

At the end of the training, we select the network having the highest empirical $W_\phi$ score, which here corresponds to choosing $\lambda=1$. We then score the initial observations $X_1, \ldots, X_{n}$ and display in Figure \ref{fig:heat_map} the points with an intensity varying from red to blue as the score increases from $0$ to $1$. The fact that the red points are on the sides of the dataset empirically validates our methodology. We represent in Fig. \ref{fig:mass_volume} the averaged mass volume curve together with standard deviation computed for $\lambda = 1$ over $B=50$ repetitions. Table \ref{restable} gathers the results averaged over $B=50$ repetitions. Notice that these results support the soundness of our approach. Indeed, the area under the $\mv$ curve is minimized and the proportion of detected outliers is high even when $n_{lowest}$ increases.

\begin{table}[h!]
\centering
\begin{tabular}{ccccc}
\hline
\vspace{0.1cm}
$n_{lowest}$ & 25  & 50 & 75 & 100  \\
\hline
\vspace{0.1cm}
$\mathrm{Acc}_{n_{lowest}}$ & $0.91 \pm 0.13$ & $0.84 \pm 0.15$   & $0.74 \pm 0.15$    &     $0.64 \pm 0.13$\\
\hline
\end{tabular}
\caption{Tabular view of the empirical accuracy $+-$ its standard deviation, when $n_{lowest}$ varies in $\{25, 50, 75, 100\}$, with $\lambda = 1$.}
\label{restable}
\end{table}

\section{Conclusion}
In this paper, we promoted a binary classification approach to the problem of learning to rank anomalies. We established a clear theoretical link between these two machine learning tasks through the study of the mass-volume curve. In particular, our procedure is robust with respect to imbalanced datasets through the choice of the parameter $p$ that is chosen initially in practice. Previous results (see \cite{CleLimVay21}) support the effectiveness of our methodology. Moreover, we illustrate our method with numerical experiments of synthetic data.

\appendix

\acks{We thank Yannick Guyonvarch for his insightful comments. Moreover, we are greatly indebted to the chair DSAIDIS of Telecom Paris and to the Région Ile-de-France for the support.}

\vskip 0.2in
\bibliography{References_ranking1,References_ranking2}
\bibliographystyle{plain}

\end{document}